\documentclass[review]{article}

\usepackage{amsfonts,amsmath,amsthm,amscd,amssymb,latexsym}
\usepackage{epsfig}
\usepackage{cancel}

\usepackage[mathlines, modulo, displaymath]{lineno} % modulo -> by default, print line numbers multiple of five

\newtheorem{theorem}{Theorem}

\newtheorem{remark}{Remark}

\newtheorem{example}{Example}

\numberwithin{equation}{section}

\usepackage{graphicx}
\usepackage{enumerate}
\usepackage{float}
\usepackage{varioref}
\usepackage{pdflscape}
\usepackage[utf8]{inputenc}
\begin{document}

\centerline{{\Large{{\textbf{Mumford representation and Riemann Roch space}}}}}

\centerline{{\Large{{\textbf{of a divisor on a hyperelliptic curve}}}}}

\bigskip
\centerline{{{Giovanni Falcone}\footnote{\textit{Address:}{Department of Mathematics and Computer Science, University of Palermo, Via Archirafi, 90123 Palermo, Italy}
\textit{email:} giovanni.falcone@unipa.it}\footnote{Supported by FFR-2023 University of Palermo.}, {Giuseppe Filippone}\footnote{\textit{Address:}{Department of Mathematics and Computer Science, University of Palermo, Via Archirafi, 90123 Palermo, Italy} \textit{email:} giuseppe.filippone01@unipa.it}\footnote{Supported by FFR-2023 University of Palermo.}}}

\begin{abstract}
For an (imaginary) hyperelliptic curve $\mathcal{H}$ of genus $g$, with a Weierstrass point $\Omega$, taken as the point at infinity, we determine a basis of the Riemann-Roch space
$\mathcal{L}(\Delta + m \Omega)$, where
$\Delta$ is of degree zero, directly from the Mumford representation of $\Delta$. This provides in turn a generating matrix of a Goppa code.
\end{abstract}

{\small{\noindent\textbf{Keywords:} Goppa codes; Riemann-Roch space;  hyperelliptic curves;  \\
 \textbf{AMS MSC codes:} 94B27 ;  14G50}}

\linenumbers

\vspace{3em}
Both the Mumford representation of a divisor $\Delta$ of degree zero on a hyperelliptic curve and the Riemann Roch space $\mathcal{L}(D)$, where $D=\Delta+m \Omega$, are the subject of a large number of papers, also due to their applications in Coding theory.

But it has not been indicated in the literature that a basis of the latter can be directly found from the former, and it is the aim of the present note to give an explicit basis of $\mathcal{L}(D)$, stressing the meaning of the Mumford representation of $\Delta$ in this context. Note that, for a nodal curve, a data structure inspired by the Mumford representation has been used for the same purpose in a recent excellent paper by Le Gluher and Spaenlehauer \cite{GluherSpaenlehauer} (the same article also contains a state of the art that we share, wishing to let this paper as short as possible).

 Using this basis, one constructs directly a generating matrix of a Goppa code over a hyperelliptic curve defined over a Galois field of characteristic $p\geq 2$. We make this for a toy model of MDS codes in Section \ref{FindH}.
Although the reduction of a divisor $D$ to its reduced Mumford form might be an inconvenient task, involving the application of the Cantor algorithm (see Remark \ref{auxiliar}), this difficulty does not occur in the construction of Goppa codes, because in that case one can directly take $D$ in the reduced form $D=\Delta+m \Omega$.

\section{Notations and reduction to the Mumford representation}
\label{sec:1}
Let $\mathsf{K}$ be the algebraic closure of the field $\mathsf{k}$ and let $\mathcal{H}$ be a hyperelliptic curve of genus $g$ over $\mathsf{k}$ with a rational Weierstrass point $\Omega$. The non-singular curve $\mathcal{H}$ is described by an affine equation of the form
\begin{equation}\label{H}
	y^2+yh(x)=f(x)
\end{equation}
where $f(x)$ is a polynomial of degree $d=2g+1$, $h(x)$ is a polynomial of degree at most $g$, and
$\Omega=[0:1:0]$ is the point at infinity of $\mathcal{H}$ (\cite{Lockhart}, Prop. 1.2). If $\operatorname{char}\mathsf{k}\ne 2$, changing
$y$ into $y-h(x)/2$, and $f(x)$ into $f(x)-h^2(x)/4$, transforms the above equation into
$$y^2=f(x),$$
whereas, if $\operatorname{char}\mathsf{k}= 2$, then it is not possible to reduce $h(x)$ to zero.

Let $D$ be a divisor of $\mathcal{H}$. Since its Riemann-Roch space  $$\mathcal{L}(D)=\{F\in \mathsf{K}(\mathcal{H}):\mathrm{div}(F)+D\mbox{ is effective}\}\cup\{0\}$$ is null both in the cases where $D$ has negative degree, and where $D$ has degree zero and
$D\not\in\mathrm{Princ}(\mathcal{H})$, whereas $\mathcal{L}(D)=\Big\langle  F_0^{-1} \Big\rangle  $ in the case where
$D=\mathrm{div}(F_0)$, from now on we will assume $D$ has positive degree $ m $, thus
$$D=P_1+P_2+\dots +P_t+(m - t)\Omega+\mathrm{div}(\psi(x,y))$$
for $t$ points $P_1,\dots , P_t$ in $\mathcal{H}$ distinct from $\Omega$, with $t\leq g$, and a suitable
$\psi(x,y)\in{\mathsf{K}}(\mathcal{H})$, that is,
any divisor class $D+\mathrm{Princ}(\mathcal{H})\in \mathrm{Div}(\mathcal{H})/\mathrm{Princ}(\mathcal{H})$ can be reduced to the form $P_1+\dots +P_t+(m - t)\Omega$.

\medskip
In order to extend the use of Mumford representation to divisors of arbitrary degree, we will apply the following:

\begin{remark}\label{auxiliar}
\emph{Note that any divisor  $$\Delta=l_1(x_1,y_1)+\dots +l_s(x_s,y_s)-(l_1+\dots +l_s)\Omega$$
	on the curve $\mathcal{H}$, of degree zero and
	such that $l_i>0$ for any index $i$, determines uniquely the polynomial
    $ a(x)= {(x-x_1)}^{l_1} \cdots {(x - x_s)}^{l_s} $
    and the polynomial $ b(x) $ which is the interpolating polynomial such that $b(x_t)=y_t$
    (hence $b^2(x)+h(x)b(x)-f(x)$ is a multiple of $a(x)$ and the degree $s-1$ of $b(x)$ is smaller than the degree of $a(x)$).
    Conversely, any pair of polynomials $a(x)$ and $b(x)$ such that $b^2(x)+h(x)b(x)-f(x)$ is a multiple
    of $a(x)$ and the degree of $b(x)$ is smaller than the degree of $a(x)$ defines such a divisor of degree zero, which is written as $\Delta=\operatorname{div}(a(x),b(x))$. Note that an intersection point of the curve with the $x$-axis is contained in the support of $\Delta$ if and only if $\operatorname{GCD}(a(x),a'(x),b(x))\ne 1$. If  $\operatorname{GCD}(a(x),a'(x),b(x))= 1$ and the degree of $a(x)$ is not greater than the genus $g$ of the curve (or equivalently, if the support of $\Delta$ contains at most $g$ points which are mutually non-opposite), one says that
    $\operatorname{div}(a(x),b(x))$ is in Mumford form (or reduced form).\\
	We remark here that any divisor $D=D_1-D_2$ (with $D_i$ effective of degree $m_i\in\mathbb{Z}$) can be written as
	$$D=\Delta+m \Omega+\operatorname{div}(\psi(x,y)),$$
	with $m = m_1-m_2$,
	for a suitable divisor $\Delta=\operatorname{div}\big(u(x),v(x)\big)$ in Mumford form, and a suitable function $\psi(x,y)$, obtained with the following argument.\\
	First, taking the vertical lines $x-x_i$ passing through the points in the support of $D_2$
    we can write $$-D_2=D_2'-2m_2\Omega-\operatorname{div}(\phi),$$
    with $\phi=\prod (x-x_i)$ and $D_2'$ effective, hence
	$$D = D_1 - D_2 = D_3-2m_2\Omega-\operatorname{div}(\phi),$$
    with $D_3=D_1+D_2'$ an effective divisor of degree $m_1+m_2$, hence of the form
    $$D_3=\operatorname{div}(a(x),b(x))+(m_1+m_2)\Omega.$$
    Secondly, applying the reduction step in Cantor's algorithm (cf. \cite{Cantor}, and \cite{koblitz} in the case where $\operatorname{char}\mathsf{k}= 2$),
    we change $D_3$ with $$D_3'=D_3-\operatorname{div}(y-b(x))=\operatorname{div}(a'(x),b' (x)),$$ which belong to the same divisor class, where
    $$a'(x)=\frac{f(x)-b(x)h(x)-b^2(x)}{a(x)}$$ and
    $$b'(x)=-h(x)-b(x)\mod a'(x).$$
    This way $\operatorname{deg} a'(x)<\operatorname{deg} a(x)$, hence after finitely many iterations one gets $\operatorname{deg} a'(x)\leq g$, and one can write
    $$D=\Delta+m \Omega+\operatorname{div}(\psi(x,y)),$$
    where $\psi(x,y)$ is the resulting function of the above reduction.\\
    	Finally, the function $$\Phi: \mathcal{L}(D)\mapsto \mathcal{L}(\Delta+m \Omega),$$ mapping $F$ onto the product
$\psi(x,y) F$, is an isomorphism. \\
Up to the latter isomorphism, we will directly assume that $D=\Delta+m \Omega$, $m > 0$.}
\end{remark}

\section{Main theorem}
\label{sec:MainThm}

In the following theorem we determine a basis of
$\mathcal{L}(D)$, with $D=\Delta + m \Omega$., and $ \Delta = \operatorname{div}(u(x),v(x)) $ is in
Mumford representation, with $t:=\operatorname{deg} u(x)\leq g$.
Also, the kind of unexpected varying, according to $m$, of its dimension becomes manifest: in order to determine $\operatorname{dim}\mathcal{L}(D)$, in \cite{boer}, Lemma 2.1 it is distinguished the case $m \geq 2g-t-1$, where $\mathrm{dim}\, \mathcal{L}(D)=m - g + 1$, and the case 	$t\leq m < 2g-t-1$, where $\mathrm{dim}\, \mathcal{L}(D)=\big\lfloor\frac{m - t}{2}\big\rfloor +1$ (cf. Remark \ref{deBoer}).

\begin{theorem}\label{Mumford} Given the hyperelliptic curve
	$\mathcal{H}$ of genus $g$ and degree $d=2g+1$ defined by \eqref{H}, given the divisor $D=\Delta+m \Omega$ of positive degree $ m $
    on $\mathcal{H}$ defined in Remark \ref{auxiliar}, with $\Delta=\operatorname{div}(u(x),v(x))$ in Mumford representation,
    let $t:=\operatorname{deg} u(x)\leq g$ and let $$\Psi(x,y)=\frac{y+ v(x)}{u(x)},$$
    for $\operatorname{char}\mathsf{k}=p>2$, and  $\Psi(x,y)=\frac{y+ v(x) + h(x)}{u(x)}$, for $p=2$.

    If $m <d-t$, then a basis of $\mathcal{L}(D)$ is provided by the set of functions $ x^i$, with $0\leq i \leq\frac{m - t}{2}$.

    If $m \geq d - t$, then a basis of $\mathcal{L}(D)$ is provided by the set of functions $ x^i$ and $\Psi(x,y)\cdot x^j$, with $	0\leq i \leq\frac{m - t}{2}$ and $0\leq j \leq\frac{m - (d - t)}{2}$.
\end{theorem}

\begin{proof}

In order to compute $\operatorname{div}(\Psi(x,y))$, recall that $\operatorname{deg} v(x)<\operatorname{deg} u(x)\leq g$ and that, in the case where $p=2$, $\operatorname{deg} h(x)\leq g$, as well.

	\medskip
	Since $ l = \operatorname{\max}(\operatorname{deg} v(x), \operatorname{deg} h(x))\leq g $,
	the degree of ${\big(- v(x) - h(x) \big)}^2$ is smaller than the degree of $ f(x) $,
	hence there are $d = 2g + 1$ intersection points of the curve $ y+v(x)+h(x)=0 $ and $\mathcal{H}$ in the affine plane, the remaining $d(l-1)$ intersection points coinciding with $\Omega$.
	More precisely, $t$ intersection points in the affine plane belong to the support of the divisor $\widehat{\Delta} =\operatorname{div} \big(u(x),w(x)\big)$ in Mumford representation, where $w(x)=-v(x)-h(x) \mod u(x)$,
	therefore
	$$\operatorname{div}\big(y+v(x)+h(x)\big)= \widehat{\Delta} + W + \big(t + d(l - 1)\big) \Omega, $$
	where
	$ W $ is the effective divisor of degree $ d - t $, whose support consists of the remaining intersection points in the affine plane. Note that, in the case $ t = 0 $,  the divisor $\Delta$ has the Mumford representation $(1,0)$ and
	the degree of $ W $ is $ d = 2g + 1 $ and the support of $ W $ coincides with the intersections
	of $\mathcal{H}$ with the curve $y+h(x)=0$.

	On the other hand, the intersection of $ u(x)=0$ and $\mathcal{H}$ is simply
	$$ \operatorname{div}\big(u(x)\big) = \Delta + \widehat{\Delta} + (t d) \Omega. $$

	Summarizing, if $t>0$, then
	\begin{equation}\label{divPsi}
		\begin{aligned}
			\operatorname{div}(\Psi(x,y)) &= \operatorname{div}\big(y+v(x)+h(x)\big) - \operatorname{div}\big(u(x)\big) =  \\
			&= W - \Delta - (d - t) \Omega
		\end{aligned}
	\end{equation}
	and, if $t=0$, then $\Delta=(1,0)$ and  $\Psi(x,y)=y+h(x)$, whence $$\mathrm{div}(\Psi(x,y))=\mathrm{div}\left(y+h(x)\right) - \mathrm{div}\left(1\right) = W - d \Omega,$$
	thus in both cases the equality (\ref{divPsi}) holds. Hence \begin{equation}\label{Psi}
	\Psi(x,y)\in\mathcal{L}(D) \mbox{ if and only if } m \geq d - t .
\end{equation}

\medskip
Let $m \geq d - t$, that is, the case where $\Psi(x,y)\in\mathcal{L}(D)$.
	First we consider the cases where either $t=0$ (hence $m \ge d = 2g+1$), or  $t=1$ (hence $m \ge d-1$), or $t \ge 2$ and $m \geq d-2$, as in these cases we know, by the theorem of Riemann-Roch, that the dimension of
	$\mathcal{L}(D)$ is $m - g + 1$. Thus, in order to prove that
	\begin{equation}\label{formula}
		\mathcal{L}(D)=\left\langle x^i,\Psi(x,y)\cdot x^j\right\rangle,\mbox{ with }\\
		0\leq i \leq\frac{m - t}{2}\mbox{ and }0\leq j \leq\frac{m - (d - t)}{2},
	\end{equation}
	it is sufficient to note that, for each of those values of the parameters $ i $ and $ j $, these functions belong to $\mathcal{L}(D)$, because
	$$1+\left\lfloor \frac{m - t}{2}\right\rfloor+1+\left\lfloor\frac{m - (d - t)}{2}  \right\rfloor=m - g + 1,$$
	and the claim will follow from dimensional reasons. Now,
	\begin{equation}\label{x^i}
	D+\mathrm{div}(x^i)=(\Delta + m \Omega) + i \cdot \mathrm{div}(x),
\end{equation}
	as well as
    \begin{equation}\label{psix^j}
        \begin{aligned}
            D+\mathrm{div}\;(\Psi(x,y)\cdot x^j) &= (\Delta + m \Omega) + j \cdot \mathrm{div}(x) + \big(W - \Delta - (d - t) \Omega\big) = \\
                                            &= W + j \cdot \mathrm{div}(x) - (d - t - m) \Omega,
        \end{aligned}
    \end{equation}
	are effective divisors, hence the functions belong to $\mathcal{L}(D)$.

	\medskip
	Secondly, we consider the case where $d-t\leq m < d - 2$. In this case, the dimension of $\mathcal{L}(D)$ is not necessarily $m - g + 1$, but still $\Psi(x,y)\in\mathcal{L}(D)$.

	If $0\leq\epsilon\leq t-2$, and if, for short, we put $m = m_\epsilon=d-2-\epsilon$, then
	$$\mathcal{L}_\epsilon:=\mathcal{L}(\Delta + m_\epsilon \Omega),$$
	hence the space
	$\mathcal{L}_0=\mathcal{L}(\Delta + (d-2)\Omega)$
	is generated, by the above case, by the functions $x^i$ and $\Psi(x,y)\cdot x^j$ with $0\leq i \leq\frac{m_0-t}{2}$ and $0\leq j \leq\frac{m_0 - (d - t)}{2}$. Of course, $\mathcal{L}_{\epsilon+1}\leq \mathcal{L}_\epsilon$, and we will see that $\operatorname{dim}(\mathcal{L}_{\epsilon+1})=\operatorname{dim}(\mathcal{L}_\epsilon)-1$.
	Indeed, by \eqref{x^i} and \eqref{psix^j},
	the functions $x^i, \Psi(x,y) x^j$ of $\mathcal{L}_\epsilon$ belong to $\mathcal{L}_{\epsilon +1}$ as long as $i \leq\frac{m_{\epsilon+1} -t}{2}$, and $j \leq\frac{m_{\epsilon+1}-(d-t)}{2}$, that is, $$\operatorname{dim}(\mathcal{L}_{\epsilon+1})=1+\left\lfloor\frac{m_{\epsilon+1} -t}{2}\right\rfloor+1+\left\lfloor\frac{m_{\epsilon+1}-(d-t)}{2}\right\rfloor,$$
	and our assertion is proved. In particular, we found that
	$\operatorname{dim}(\mathcal{L}_{\epsilon+1})=\operatorname{dim}(\mathcal{L}_\epsilon)-1$, because $m_{\epsilon+1}=m_\epsilon-1$ and
	$$\operatorname{dim}(\mathcal{L}_{\epsilon})=1+\left\lfloor\frac{m_{\epsilon} -t}{2}\right\rfloor+1+\left\lfloor\frac{m_{\epsilon}-(d-t)}{2}\right\rfloor,$$ where the missing function is, once for one, $x^i$ or $\Psi(x,y) x^j$, because $d$ is odd and changes the parity of ${m_{\epsilon+1} -t}$ in that of $m_{\epsilon+1}-(d-t)$.

	\bigskip\noindent
	Now we consider the cases where $m < d-t$, that is, the cases where, by \eqref{H}, $\Psi(x,y)\not\in\mathcal{L}(D)$. If $t=0$ and $m\in\{d-2, d-1\}$, or if $t=1$ and $m=d-2$, then on the one hand $\lfloor\frac{m - t}{2}\rfloor=m-g$ and, on the other hand,  by the theorem of Riemann-Roch, the dimension of $\mathcal{L}(D)$ is $m - g + 1$. Thus, by dimensional reason,
	$\mathcal{L}(D)=\left\langle x^i\right\rangle$, where $0\leq i \leq\lfloor\frac{m - t}{2}\rfloor$.

	In order to prove that
	$\mathcal{L}(D)=\left\langle x^i\right\rangle$, where $0\leq i \leq\frac{(m - t)}{2}$ also in the remaining cases where either $t=0, 1$ and $m <d-2$, or $2\leq t\leq m <d-t$, write
	$m = m_\epsilon=d-t-\epsilon$ with $1\leq\epsilon\leq d-2t$, and again put, for short,
	$$\mathcal{L}_\epsilon:=\mathcal{L}(\Delta + m_\epsilon \Omega).$$
	Note that appending the value $\epsilon=0$, that is, considering also the case where $m = m_0=d-t$, by \eqref{formula} we have
	$\mathcal{L}_0=\left\langle x^i,\Psi(x,y)\right\rangle$,
	with $0\leq i \leq\frac{m_0-t}{2}$.

	Of course, $\mathcal{L}_{\epsilon+1}\leq \mathcal{L}_\epsilon$ for any $0\leq\epsilon\leq d-2t$, but in this case we will see that \begin{equation}\label{dim}
\operatorname{dim}(\mathcal{L}_{\epsilon+1})=
\left\{\begin{array}{ll}
	\operatorname{dim}(\mathcal{L}_\epsilon)&\mbox{if $m_\epsilon - t$ is odd, }\\
	\operatorname{dim}(\mathcal{L}_\epsilon)-1&\mbox{if $m_\epsilon - t$ is even. }\end{array}	\right.
	\end{equation}
	 Indeed, by \eqref{Psi} $\Psi(x,y)\not\in\mathcal{L}_\epsilon$ as soon as $\epsilon>0$,
and since, by \eqref{x^i},
the functions $x^i$ of $\mathcal{L}_\epsilon$ belong to $\mathcal{L}_{\epsilon +1}$ as long as $i \leq\frac{m_{\epsilon+1} -t}{2}$, we see that $$\operatorname{dim}(\mathcal{L}_{\epsilon+1})=1+\left\lfloor\frac{m_{\epsilon+1} -t}{2}\right\rfloor,$$
and we get the equalities in \eqref{dim},
because $m_{\epsilon+1}=m_\epsilon-1$. But this equality shows, as well, that the theorem is true for any value of $ m $.
\end{proof}

\begin{remark}\label{deBoer}
	\emph{It is remarkable that the bounds in \cite{boer}, Lemma 2.1 are different from the ones above: for $m = 2g-t, 2g-t-1$, in our theorem we find $\operatorname{dim}\mathcal{L}(D)=1+\left\lfloor\frac{m - t}{2}\right\rfloor$, whereas in  \cite{boer}, Lemma 2.1. we read $\operatorname{dim}\mathcal{L}(D)=m - g + 1$. Of course, the two values coincides exactly for  $m = 2g-t, 2g-t-1$. \\
	\indent	In particular, the necessary condition in  \cite{boer}, Lemma 2.1. to have $\operatorname{dim}\mathcal{L}(D)\ne m - g + 1$, that is, $m <d-t-2$, is also sufficient.\\
%i) Recall that $\Omega$ was a given Weierstrass point of the curve $\mathcal{H}$, and in fact
%		$\mathcal{L}\big((2g-2\epsilon)\Omega\big)=\mathcal{L}\big((2(g-\epsilon)+1)\Omega\big)$ have both dimension
%		$g-\epsilon+1$. The sequence $\mathrm{dim}\,\mathcal{L}(m \Omega)$, where $n \ge 0$, is therefore:
%		$$1,1,2,2,3,3,\dots, g-1,g-1, g, g, g+1,g+2,g+3\dots$$
%		The numerical semigroup of non-gaps is therefore that of the natural numbers without the odd numbers smaller than $2g$.\\ ii)
	\indent	An interesting phenomenon occurs when $g<m <2g-1$ and $t\in\{g,g-1,g-2\}$, because in these cases $m \geq d-t-2$, hence $\mathrm{dim}\,\mathcal{L}(D)=m - g + 1$, regardless of the theorem of Riemann-Roch.}
\end{remark}

\section{Applications to coding theory}\label{FindH}
In this section we assume $ \mathsf{k} = \operatorname{GF}(p^c) $, where $ p \ge 2 $ is a prime number and $ c $ a positive integer.

Note that, for any polynomials $u(x)$ and $v(x)$, with $u(x)$ of degree $t$, and $v(x)$ of degree smaller than $t$ (and, if $p=2$, for any arbitrary non-zero polynomial $h(x)$), there is a hyperelliptic curve of arbitrary genus $g\geq \max\{t, \operatorname{deg}(h)\}$, of equation $y^2=v^2(x)+v(x)h(x) - c(x)u(x)$, for each polynomial $c(x)$ of degree $2g+1-t$, passing through the support of $D=\Delta + m \Omega$, with $\Delta=\operatorname{div}\big(u(x),v(x)\big)$ in Mumford representation, and all of these curves determine the same Riemann Roch space $\mathcal{L}(D)$ for $D$. That is, in order to give the basis of the space $\mathcal{L}(D)$ one does not have to know the curve containing the support of $D$. Note also that one does not need to give explicitly the points in the support of $D$, a sensible advantage in the construction of AG-Goppa codes, as we will see in Example \ref{example1}.
In that Example, we compute the generating matrix of a toy model of a Goppa code of length $n=10$ and dimension $k = 5$, arising from a hyperelliptic curves of genus $g=11$, and which is a MDS code, although here the Goppa lower bound is equal to $-5$.

\begin{remark}
	\emph{	Note that, for $p\leq\frac{m - t}{2}$, the polynomials $x$ and $x^{p^c}$ in the basis of $\mathcal{L}(D)$ take the same values in the field $\mathsf{k}=\mathrm{GF}(p^c)$, and the same occurs, for $p\leq\frac{m - (d - t)}{2}$, to the polynomials $\Psi(x,y) x$ and $\Psi(x,y) x^{p^c}$. This fact must be taken into account, for instance, when constructing a Goppa code. }
\end{remark}

\begin{theorem} Let $ \mathsf{k} $ be a field of characteristic $ p \ge 2 $, let
	 $u(x)$ be a monic polynomial of degree $t$ and $v(x)$ be a polynomial with $\operatorname{deg}(v)<t$,
    such that $ \operatorname{GCD}(u(x),u'(x),v(x)) = 1 $, and let $P_s=(x_s,y_s)$ be $n$ pairs such that $u(x_s)\ne 0$,
    for any $s=1,\dots ,n$.

    If $ g \ge t $, then, for any $g-t+2 \le k < n $ the matrix $G=(\gamma_{rs})$
	\begin{equation}\left\{
		\begin{array}{ll}
			\gamma_{rs} = x_s^{r-1}&\mbox{ for }1\leq r\leq \eta+1\\
			&\\
			\gamma_{rs} =  \Psi(x_s,y_s)\cdot x_s^{r-\eta} &\mbox{ for }\eta+2\leq r\leq k
		\end{array} \right.\ \left(\mbox{where }\eta=\left\lfloor\frac{k + g - 1 - t}{2}\right\rfloor\right)
	\end{equation}
	is the generating matrix of a $[n, k, \delta]$ Goppa code, with $ n - k + 1 - g \leq \delta \leq n - k + 1 $, and with $\Psi(x_s,y_s)=\frac{y_s + v(x_s) + h(x_s)}{u(x_s)}$, where $h(x)=0$, if $p>2$, or $h(x)$ is an arbitrary non-zero polynomial with $\operatorname{deg}(h)\leq g$, if $p=2$.
\end{theorem}

\begin{proof}
    Let $ c(x) $ be a polynomial of degree $ 2g + 1 - t $ such that
    $$ c(x_s) = \frac{{v(x_s)}^2 + h(x_s) v(x_s) - y_s^2 - y_s h(x_s)}{u(x_s)}, $$
    for any $ (x_s, y_s) $ with $s=1,\dots ,n$.

	Hence, there is an hyperelliptic curve of genus $ g $ of equation
	$$ y^2 + y h(x) = f(x) = {v(x)}^2 + h(x) v(x) - c(x) u(x), $$
	passing through the $ n $ points $ (x_s, y_s) $ and the points belonging to the support of the divisor
    $ \operatorname{div}(u(x), v(x)) $.

    The claims follows from the fact that the functions taken into account in the theorem give in turn a basis
    of the Riemann-Roch space $ \mathcal{L}(D) $, where $ D = \operatorname{div}(u(x), v(x)) + (k + g - 1) \Omega $, whose
    dimension is $ k $.
\end{proof}

\begin{remark}
	\emph{	Note that, as long as $ k < g - t + 2 $ and the $n$ points $P_s=(x_s, y_s)$ where we evaluate the functions of the basis of $\mathcal{L}(D)$ have different {absciss\ae} $x_s$, the Goppa code coincides with the $[n,k,n-k+1]$ Reed-Solomon code on the $n$ values $\{x_1,\dots ,x_n\}\subset\mathsf{k}$.}
\end{remark}

%\begin{example} \emph{On the hyperelliptic curve
%		$\mathcal{H}: y^2=x^5+13x^4+5x^3+11x^2+5x+15$ over the field $\operatorname{GF}(17)$, we choose the divisior
%		$D=(8,0)+3\Omega=(x-8,0)+4\Omega$. Hence one has $\Psi(x,y)=\frac{y}{x-8}$. From Theorem \ref{Mumford} we obtain $\mathcal{L}(D)=\left\langle 1, x, \Psi(x,y)\right\rangle $ since $m = 4$.
%		With respect to the divisor $G=(0,7)+(0,10)+(1,4)+(1,13)+(3,8)+(3,9)+(5,1)+(5,16)+(9,1)+(9,16)+ (15,7)+(15,10)\in \mathcal{H}$, not in the support of $D$, the generating matrix of the
%		$(12,3,10)_{17}$-MDS code is
%		$\mathcal{C}_{\mathcal{L}}(D,G)=\left( \begin{array}{cccccccccccc}
%			1&1&1&1&1&1&1&1&1&1&1&1\\
%			0&0&1&1&3&3&5&5&9&9&15&15\\
%			14&3&14&3&12&5&11&6&1&16&1&16 \end{array} \right)$.} \end{example}

\begin{example}\label{example1} \emph{Let $\mathsf{k}=\operatorname{GF}(101)$, choose a pair of polynomials $(u(x),v(x))$ with $\operatorname{GCD}(u(x),u'(x),v(x))=1$, for instance $(u(x),v(x))= (x^{11} + 1, x^6 + 1)$, and consider the function  $$\Psi(x,y)=\frac{y+v(x)}{u(x)}=\frac{y+x^6+1}{x^{11}+1}.$$
%		$D=\operatorname{div}(u(x), v(x))+4\Omega = (x^{11} + 1, x^6 + 1)+4\Omega$.
 %Hence one has $\Psi(x,y)=\frac{y+x^6+1}{x^{11}+1}$.
 Choose five pairs $(x_s, y_s)$ such that $x_r\ne x_l$ whenever $r\ne l$, such that $u(x_s)\ne 0$ for any index $s$, for instance
		$ (15, 45)$,  $ (53, 48)$,  $ (58, 10)$, $ (64, 13)$,  $ (80, 2)$.  Evaluating the functions $\left\{ 1, x, x^2, \Psi(x,y),x\Psi(x,y)\right\} $ on the ten points 	$ (15, \pm 45)$,  $ (53, \pm 48)$,  $ (80, \pm 2)$,  $ (58, \pm 10)$, $ (64, \pm 13)$, one obtains a matrix
		$$G:= \left( \begin{array}{rrrrrrrrrr}
			1 & 1 & 1 & 1 & 1 & 1 & 1 & 1 & 1 & 1 \\
			15 & 15 & 53 & 53 & 80 & 80 & 58 & 58 & 64 & 64 \\
			23 & 23 & 82 & 82 & 37 & 37 & 31 & 31 & 56 & 56 \\
			73 & 41 & 35 & 92 & 1 & 45 & 99 & 71 & 48 & 21 \\
			85 & 9 & 37 & 28 & 80 & 65 & 86 & 78 & 42 & 31
		\end{array} \right)$$ 	which is a generating matrix
	of a code $\mathcal{C}$ of lenght $10$ and dimension $5$. A direct computation would show that the ten chosen points are such that all the $5\times 5$ minors of $G$ have full rank, hence the minimal distance is $6$, that is, $\mathcal{C}$ is a $[10,5,6]_{101}$ MDS code. \\
	In order to give the equation of a hyperelliptic curve  $\mathcal{H}$ realizing the above code as an AG-Goppa code, defined by $D=\operatorname{div}(u(x),v(x))+15\Omega$ by evaluating the functions in $\mathcal{L}(D)$ on the above five points $(x_s, y_s)$, we note that the genus $g$ of $\mathcal{H}$ must be equal at least to the degree of $u(x)$. With $g$ equal to the degree of $u(x)$, hence with the degree of $\mathcal{H}$ equal to $23$, we need eight further points, because $\mathcal{H}$ passes through the five points $(x_s, y_s)$ and through the eleven points (in the affine plane) of the support of $\operatorname{div}(u(x),v(x))$. Choose arbitrarily eight pairs $(x_s, y_s)$ (now with $s=6,\dots ,13$) such that $ u(x_s) \ne 0 $, for instance $(48, 80)$, $(58, 91)$, $(64, 88)$, $(89, 16)$, $(95, 33)$, $(53, 4)$, $(51, 85)$, $(71, 35)$.\\
	With this choice, the curve $\mathcal{H}$ defined by the equation $$y^2=v^2(x)-c(x)u(x),$$
	where $c(x)$ is the polynomial such that  $$c(x_s)=\frac{{v(x_s)}^2 - y_s^2}{u(x_s)},$$
	for $s=1,\dots, 13$, has degree $23$, passes through the $13$ points $(x_s, y_s)$ and the eleven points (in the affine plane) of the support of $\operatorname{div}(u(x),v(x))$, thus realizing the $[10,5,6]_{101}$ code as the AG Goppa code defined by $\mathcal{L}(D)$ and the ten points $(x_s,\pm y_s)$, for $s=1,\dots ,5$.} \end{example}


\begin{thebibliography}{}


%\bibitem{BLP} D.J. Bernstein, T. Lange, C. Peters, Attacking and Defending the McEliece Cryptosystem, in: J. Buchmann, J. Ding  (Eds.), Post-Quantum Cryptography, PQCrypto 2008, LNCS, vol. 5299.,  Springer, Berlin, Heidelberg, 2008, pp. 31–46.

%\bibitem{boer0} M.A. de Boer, MDS Codes from Hyperelliptic Curves. in: R. Pellikaan, M. Perret, S.G. Vl{\"a}dut  (Eds.), Arithmetic, Geometry and Coding Theory, Walter de Gruyter, Berlin, New York, 1996, pp. 23-34.

%\bibitem{boer1} M.A. de Boer, Codes: their parameters and geometry, Eindhoven, Technische Universiteit Eindhoven, 1997, DOI: 10.6100/IR492527.

\bibitem{boer} M.A. de Boer, The generalized Hamming weights of some hyperelliptic codes, J. Pure Appl. Algebra 123 (1998) 153-163.

%\bibitem{BrillNoether} A. von Brill,  M. Noether, \"Uber die algebraischen Functionen und ihre Anwendung in der Geometrie, Math. Annalen 7 (1874) 269-316.


\bibitem{Cantor} D.G. Cantor, Computing in the Jacobian of a Hyperelliptic Curve,
Math. Comp. 48 (1987) 95-101.

%\bibitem{CastellanosFanali} A.S. Castellanos, G.C. Tizziotti, Two-Point AG Codes on the GK Maximal Curves, IEEE Trans. Inf. Theory 62 (2016) 681-686.

%\bibitem{DinhMooreRussell} H. Dinh, C. Moore, A. Russell, McEliece and Niederreiter Cryptosystems That Resist Quantum Fourier Sampling Attacks, in: P. Rogaway (Ed.), Advances in Cryptology – CRYPTO 2011, CRYPTO 2011, LNCS, vol. 6841., Springer, Berlin, Heidelberg, 2011, pp. 761-779.

%\bibitem{FanaliGiulietti} S. Fanali, M. Giulietti, One-Point AG Codes on the GK Maximal Curves, IEEE Trans. Inf. Theory 56 (2010) 202-210.


%\bibitem{G2} V.D. Goppa, Algebraic-geometric codes, Izv. Akad. Nauk SSSR Ser. Mat. 46 (1982)
% 762-781. (in Russian)

%\bibitem{HFF} G. Falcone, \'A. Figula and C. Hannusch, On the generating matrices of Goppa codes over hyperelliptic curves, J. Ramanutan Math. Soc. 37 (3), 273–279 (2022).


%\bibitem{Hess}  F. Hess, Computing Riemann-Roch Spaces in Algebraic Function Fields and Related Topics, J. Symbolic Comp. 33 (2002) 425--445.

%\bibitem{HuangIerardi} M. Huang, D. Ierardi, Efficient algorithms for the Riemann-Roch problem and for addition in the Jacobian of a curve, J. Symbolic Comp. 18 (1994) 519-539.

%\bibitem{JanwaMoreno} H. Janwa, O. Moreno, McEliece public key cryptosystems using algebraic-geometric codes, Des. Codes Crypt. 8  (1996) 293–307.

\bibitem{koblitz} N. Koblitz, Hyperelliptic Cryptosystems,  J. Cryptology 1 (1989) 139–150.

%\bibitem{KorchmarosNagyTimpanella} G. Korchm\'aros, G.P. Nagy, M. Timpanella, Codes and Gap Sequences of Hermitian Curves,  IEEE Trans. Inf. Theory 66 (2020) 3547--3554.

%\bibitem{KorchmarosSpeziali} G. Korchm\'aros, P. Speziali, Hermitian codes with automorphism group isomorphic to $\mathrm{PGL}(2, q)$ with $q$ odd, Finite Fields Their Appl. 44 (2017) 1-17.

%\bibitem{kuroki} J. Kuroki, M. Gonda, K. Matsuo, J. Chao, S. Tsutii, Fast Genus Three Hyperelliptic Curve Cryptosystems, In The 2002 Symposium on Cryptography and Information Security, Japan - SCIS 2002, 2002.

%\bibitem{Lange2}  T. Lange, Formulae for Arithmetic on Genus 2 Hyperelliptic Curves, AAECC 15 (2005) 295-328.

%\bibitem{brigand} D. Le Brigand, Decoding of codes on hyperelliptic curves, in: G. Cohen,  P. Charpin (Eds.), EUROCODE '90, EUROCODE 1990, LNCS, vol. 514., Springer, Berlin, Heidelberg, 1990, pp. 126--134.

\bibitem{GluherSpaenlehauer} A. Le Gluher,  P.-J. Spaenlehauer, A fast randomized geometric algorithm for computing Riemann-Roch spaces, Math. Comp. 89 (2020) 2399-2433.

\bibitem{Lockhart} P. Lockhart, On the discriminant of a hyperelliptic curve, Trans. Amer. Math. Soc. 342 (1994) 729-752.

%\bibitem{MMPR} I. Márquez-Corbella, E. Martínez-Moro, R. Pellikaan, D. Ruano, Computational aspects of retrieving a representation of an algebraic geometry code, J. Symbolic Comp. 64 (2014) 67-87.

%\bibitem{matsuo} K. Matsuo, J. Chao, S. Tsutii, Fast Genus Two Hyperelliptic Curve Cryptosystems, ISEC 2001-31, %IEICE 2001.

%\bibitem{Mc} R.J. McEliece, A Public-Key Cryptosystem Based On Algebraic Coding Theory, DSN Progress Report 44 (1978) 114-116.

%\bibitem{niehage2} A. Niehage, Nonbinary Quantum Goppa Codes Exceeding the Quantum Gilbert-Varshamov Bound, Quantum Inf. Process 6 (2007) 143-158.

%\bibitem{pelzl} J. Pelzl, T. Wollinger, J. Guatardo, C. Paar, Hyperelliptic curves cryptosystems: closing the performance gap to elliptic curves,  in: C.D. Walter, Ç.K. Koç, C. Paar  (Eds.) Cryptographic Hardware and Embedded Systems - CHES 2003, CHES 2003, LNCS, vol. 2779., Springer, Berlin, Heidelberg, 2003, pp. 351-365.


%\bibitem{stepanov} S.A. Stepanov, Codes on fibre products of hyperelliptic curves,
%Disc. Math. Appl. 7 (1997) 77--88.

%\bibitem{Stichtenoth} H. Stichtenoth, Algebraic function fields and codes, Springer, Berlin, Heidelberg, 2009.

%\bibitem{sutherland} A.V. Sutherland, Fast Jacobian arithmetic for hyperelliptic curves of genus $3$, Thirteenth Algorithmic Number Theory Symposium (ANTS XIII), Open Book Series 2, 2019, pp. 425-442.

%\bibitem{yaghoobian} T. Yaghoobian, I.F. Blake, Codes from hyperelliptic curves, in: Proc. 30-th Allerton Conf. %Comm., Control, and Computing, Monticello, IL., October 1992.




\end{thebibliography}
\end{document}